\documentclass[11pt]{article}

\usepackage{amssymb,latexsym}
\usepackage{amsmath}
\usepackage{cite}
\usepackage{eucal}

\title{Morita equivalence of Fedosov star products and deformed
  Hermitian vector bundles} 

\author{\textbf{Stefan Waldmann\thanks{Stefan.Waldmann@physik.uni-freiburg.de}}  
  \\[0.5cm]
  Fakult{\"a}t f{\"u}r Physik\\ 
  Albert-Ludwigs-Universit{\"a}t Freiburg\\
  Hermann Herder Stra{\ss}e 3\\
  D 79104 Freiburg\\
  Germany
  }

\date{January 2002\\[0.5cm] FR-THEP 2002/1}

\newcommand{\im} {{\mathrm i}}
\newcommand{\eu} {{\mathrm e}}

\newcommand{\cc}[1]      {\overline{{#1}}}

\newcommand{\id}         {{\mathsf{id}}}

\newcommand{\ad}         {{\mathrm{ad}}}

\newcommand{\End}        {{\mathsf{End}}}

\newcommand{\Weyl}       {{\mathcal{W}}}
\newcommand{\WL}         {{\mathcal{W}\!\otimes\!\Lambda^{\!\bullet}}}
\newcommand{\WLE}        {{\mathcal{W}\!\otimes\!\Lambda^{\!\bullet}\!\otimes\!\mathcal{E}}}
\newcommand{\WLEnd}      {{\mathcal{W}\!\otimes\!\Lambda^{\!\bullet}\!\otimes\!\mathcal{E}nd(\mathcal{E})}}
\newcommand{\CWeyl}      {{\hat{\mathcal{W}}}}
\newcommand{\CWL}        {{\hat{\mathcal{W}}\!\otimes\!\Lambda^{\!\bullet}}}
\newcommand{\CWLE}       {{\hat{\mathcal{W}}\!\otimes\!\Lambda^{\!\bullet}\!\otimes\!\mathcal{E}}}
\newcommand{\CWLEnd}     {{\hat{\mathcal{W}}\!\otimes\!\Lambda^{\!\bullet}\!\otimes\!\mathcal{E}nd(\mathcal{E})}}
\newcommand{\WLO}        {\mathcal{W}\!\otimes\!\Lambda}

\newcommand{\womega}     {\widetilde{\omega}}

\newcommand{\nablaE}     {\nabla^E}
\newcommand{\nablaEnd}   {\nabla^{\End(E)}}

\newcommand{\degs}       {\mathrm{deg}_{\mathrm{s}}}
\newcommand{\dega}       {\mathrm{deg}_{\mathrm{a}}}
\newcommand{\degl}       {\mathrm{deg}_{\lambda}}
\newcommand{\Deg}        {\mathrm{Deg}}

\newtheorem{lemma}{Lemma}

\newtheorem{theorem}{Theorem}
\newtheorem{corollary}{Corollary}

\newtheorem{remark}{Remark}

\newenvironment{proof}{\small{\sc Proof:}}{{\hspace*{\fill} $\square$\\}}

\numberwithin{equation}{section}

\allowdisplaybreaks

\begin{document}

\maketitle

\begin{abstract}
    Based on the usual Fedosov construction of star products for a
    symplectic manifold $M$ we give a simple geometric construction of
    a bimodule deformation for the sections of a vector bundle over
    $M$ starting with a symplectic connection on $M$ and a connection
    for $E$. In the case of a line bundle this gives a Morita
    equivalence bimodule where the relation between the characteristic
    classes of the Morita equivalent star products can be found very
    easily in this framework. Moreover, we also discuss the case of a
    Hermitian vector bundle and give a Fedosov construction of the
    deformation of the Hermitian fiber metric.
\end{abstract}


\newpage

%
%

\section{Introduction}

Deformation quantization as introduced in \cite{bayen.et.al:1978} has
proved to be an extremely successful framework for the problem of
quantization: the existence of the associative deformation of the
classical observable algebra, the `star product', is well established
for the case of a symplectic phase space
\cite{dewilde.lecomte:1983b,fedosov:1996,omori.maeda.yoshioka:1991} as
well as for the more general Poisson case
\cite{kontsevich:1997:pre}. Moreover, star products have been
classified up to equivalence
\cite{nest.tsygan:1995a,bertelson.cahen.gutt:1997,kontsevich:1997:pre}
in terms of geometrical data on the phase space. For several physical
applications one also needs to represent the deformed observable
algebra on a pre-Hilbert space. This led to the development of a
representation theory for star products starting with
\cite{bordemann.waldmann:1998}. Recent reviews as well as further
references may be found in
\cite{gutt:2000,waldmann:2001b:pre,weinstein:1994,dito.sternheimer:2001:pre}.

Having a reasonable notion for a representation theory, a natural
question is whether two star product algebras have the `same'
representation theory. This question was made precise in
\cite{bursztyn.waldmann:2001a} using a notion of Morita equivalence
very similar to and in fact generalizing Rieffel's notion of strong
Morita equivalence for $C^*$-algebras \cite{rieffel:1974b}. The
classification of star products up to Morita equivalence was achieved
in \cite{bursztyn.waldmann:2001:pre} for the symplectic case. In
the particular case of cotangent bundles it leads to a physical
interpretation of Morita equivalence as Dirac's quantization condition
for magnetic charges. Beside this more `conservative' occurrence of
Morita equivalence in deformation quantization, Morita equivalence of
star products also appears in non-commutative gauge theories, see
\cite{jurco.schupp.wess:2001b:pre}.

The purpose of this paper is to give an alternative and more geometric
construction of the deformation of vector bundles $E \to M$ as
introduced in \cite{bursztyn.waldmann:2000b} which are the basic
ingredients for Morita equivalence of star products. Here we shall use
a `Fedosov-like' construction of the bimodule structure and give
thereby a simple description of Morita equivalence for Fedosov star
products.

The paper is organized as follows: first we recall the basic
structures needed for Fedosov's approach to star products in Section~2
and 3, where we also introduce the fiberwise bimodule structure. In
the next two sections we show how the usual Fedosov derivatives, which
lead to deformations of $C^\infty(M)$ and $\Gamma^\infty(\End(E))$,
can be used to obtain also a Fedosov derivative for the vector bundle
itself. This will allow to define the deformed bimodule structure as
well as an easy identification of the characteristic classes of the
involved star products. In Section~6 we demonstrate how a Hermitian
fiber metric can be deformed in this framework and Section~7 contains
a conclusion with some further questions arising in this context.

\medskip
\noindent
\textbf{Acknowledgements:}
I would like to thank Henrique Bursztyn and Nikolai Neumaier for a
careful reading of the manuscript and many comments and useful
discussions. Moreover, I would like to thank the participants of the
Warwick workshop on topology, operads and quantisation for their
remarks and encouraging comments.

%
%

\section{Preliminaries on the Fedosov construction}
\label{sec:prelim}

The aim of this section is to recall the basics of Fedosov's
construction and to set up our notation, where we mainly follow
\cite{bordemann.waldmann:1997a}. In the following $(M,\omega)$ 
is a symplectic manifold, $\nabla$ a symplectic torsion-free
connection, $E \to M$ a complex vector bundle, and $\nablaE$ a
connection for $E$. By $\End(E) \to M$ we denote the endomorphism
bundle of $E$ and $\nablaEnd$ is the induced connection for
$\End(E)$ coming from $\nablaE$. The starting point for the Fedosov
construction are the following $\mathbb{C}[[\lambda]]$-modules:
\begin{gather}
    \label{eq:Wdef}
    \Weyl 
    := \prod_{s=0}^\infty 
    \Gamma^\infty (\mbox{$\bigvee$}^s T^*M) [[\lambda]]
    \\
    \label{eq:WLdef}
    \WL 
    := \prod_{s=0}^\infty \Gamma^\infty 
    (\mbox{$\bigvee$}^s T^*M \otimes \mbox{$\bigwedge$}^{\!\bullet}\, T^*M)
    [[\lambda]]
    \\
    \label{eq:WLEdef}
    \WLE 
    := \prod_{s=0}^\infty \Gamma^\infty
    (\mbox{$\bigvee$}^s T^*M \otimes 
    \mbox{$\bigwedge$}^{\!\bullet}\, T^*M \otimes E) [[\lambda]]
    \\
    \label{eq:WLEnddef}
    \WLEnd
    := \prod_{s=0}^\infty \Gamma^\infty 
    (\mbox{$\bigvee$}^s T^*M 
    \otimes \mbox{$\bigwedge$}^{\!\bullet}\, T^*M \otimes \End(E)) [[\lambda]]
\end{gather}
Clearly, $\WL$ is a super-commutative associative algebra for the
(anti-)symmetric tensor product and $\WLEnd$ is associative but
non-commutative unless $E$ is a line bundle. Moreover, $\WLE$ is a
$\WLEnd$ left 
module and a $\WL$ right module in the obvious way and both module
actions commute. Thus it is a bimodule. The \emph{degree derivations}
$\degs$, $\dega$, $\degl$ and the \emph{total degree} 
$\Deg = \degl + 2\degs$ are defined in 
the usual way and yield (module-) derivations for the (yet still
undeformed) fiberwise products. Finally we have the operators 
$\delta:= (1 \otimes dx^i)i_s(\partial_i)$ and 
$\delta^* := (dx^i \otimes 1) i_a (\partial_i)$, acting on $\WL$,
$\WLE$ and $\WLEnd$, where we write $(1 \otimes dx^i)$ to emphasize that
the one-form $dx^i$ is considered to be `anti-symmetric' while
$(dx^i \otimes 1)$ refers to a `symmetric' one-form. The maps
$i_s(\partial_i)$ and $i_a(\partial_i)$ are the symmetric and
anti-symmetric insertion maps, respectively. One has 
$\delta^2 = 0 = (\delta^*)^2$ and 
$\delta \delta^* + \delta^* \delta = \degs + \dega$. Defining the
symbol map $\sigma$ as projection onto the part of symmetric and
antisymmetric degree $0$  one has
\begin{equation}
    \label{eq:hodge}
    \delta \delta^{-1} + \delta^{-1} \delta + \sigma = \id,
\end{equation}
where $\delta^{-1} a = \frac{1}{k+l}\, \delta^*a$ for $\degs a = ka$,
$\dega a = la$ with $k+l \ne 0$ and $\delta^{-1} a = 0$ else.

We shall need a slight generalization of $\Weyl$ in the following. Let
$\CWeyl$ denote the space of formal Laurent series in the total degree
$\Deg$ such that $a \in \CWeyl$ is of the form 
$a = \sum_{r=N}^\infty a^{(r)}$ with $N \in \mathbb{Z}$ and each
$a^{(r)}$ may contain arbitrarily high negative powers of $\lambda$ as
long as they are compensated by the symmetric degree such that still
$\Deg a^{(r)} = r a^{(r)}$. It is clear that all the algebraic
operations are still defined on $\CWeyl$, $\CWL$, $\CWLE$, and
$\CWLEnd$, respectively, and obey the same algebraic identities. Note
that now $\sigma$ maps $\CWL$ onto the space
$C^\infty(M)(\!(\lambda)\!)$ of formal Laurent series in $\lambda$.

The connections $\nabla$ and $\nablaE$ extend to super-derivations of
anti-symmetric degree $+1$
\begin{equation}
    \label{eq:DDef}
    D: \WL \to \WL^{+1}
\end{equation}
\begin{equation}
    \label{eq:DEDef}
    D^E: \WLE \to
    \mathcal{W}\!\otimes\!\mbox{$\bigwedge$}^{\!\bullet+1}\!\otimes\!\mathcal{E}
\end{equation}
\begin{equation}
    \label{eq:DEnd}
    D': \WLEnd \to
    \mathcal{W}\!\otimes\!\mbox{$\bigwedge$}^{\!\bullet+1}\!\otimes\!\mathcal{E}nd(\mathcal{E})
\end{equation}
in the usual way. In particular we have the compatibility conditions
\begin{equation}
    \label{eq:DDEDEndComp}
    D^E (a\Psi) = (D'a)\Psi + (-1)^{\dega a} a (D^E\Psi)
    \quad
    \textrm{and}
    \quad
    D^E(\Psi b) = (D^E\Psi) b + (-1)^{\dega \Psi} \Psi (Db)
\end{equation}
for $a \in \WLEnd$, $\Psi \in \WLE$, and $b \in \WL$. A straightforward
computation shows that
\begin{equation}
    \label{eq:deltaDCom}
    [\delta, D] = 0,
    \quad
    [\delta, D^E] = 0,
    \quad
    \textrm{and}
    \quad
    [\delta, D'] = 0.
\end{equation}
The symplectic form $\omega$ can be viewed as element in $\WLO^2$ but
we can also view it as $\womega \in \WLO^1$ of symmetric degree $1$
and antisymmetric degree $1$. Then we obtain
\begin{equation}
    \label{eq:deltaomega}
    \delta \womega = 2 \omega, 
    \quad
    \delta^{-1}\womega = 0, 
    \quad
    \textrm{and}
    \quad
    D\womega = 0,
\end{equation}
since $\nabla$ is a symplectic connection. Next we compute the squares
of the covariant derivatives where we obtain curvature
contributions. We define $R \in \WLO^2$ of symmetric degree $2$ by
\begin{equation}
    \label{eq:RDef}
    R = 
    \frac{1}{4} \, \omega_{lm} R^m_{kij} 
    \, dx^l \vee dx^k \otimes dx^i \wedge dx^j
\end{equation}
and $R^E \in \WLO^2\!\otimes\!\mathcal{E}nd(\mathcal{E})$ of symmetric
degree $0$ by
\begin{equation}
    \label{eq:REDef}
    R^E = \frac{1}{2} \, dx^i \wedge dx^j \otimes R^E_{ij},
\end{equation}
where $R^m_{kij}$ are the components of the curvature of $\nabla$ and
$R^E_{ij} = R^E(\partial_i,\partial_j)$ are the curvature
endomorphisms of $\nablaE$. Clearly the above local formulas define
global objects. A straightforward computation gives the following
lemma:
\begin{lemma}
    \label{lemma:Dsquares}
    One has the following relations:
    \begin{equation}
        \label{eq:Dsquare}
        D^2 = 
        - \frac{1}{2} \, R^m_{kil} 
        \, dx^k \otimes dx^i \wedge dx^j 
        \; i_s(\partial_m)
    \end{equation}
    \begin{equation}
        \label{eq:DEsquare}
        (D^E)^2 = 
        - \frac{1}{2} \, R^m_{kil} 
        \, dx^k \otimes dx^i \wedge dx^j 
        \; i_s(\partial_m) + R^E
    \end{equation}
    \begin{equation}
        \label{eq:DEndsquare}
        (D')^2 = 
        - \frac{1}{2} \, R^m_{kil} 
        \, dx^k \otimes dx^i \wedge dx^j 
        \; i_s(\partial_m) + [R^E, \, \cdot\,]
    \end{equation}
    The curvature elements fulfill the Bianchi identities
    \begin{equation}
        \label{eq:bianchi}
        \delta R = 0,
        \quad
        \delta R^E = 0
        \quad
        \textrm{and}
        \quad
        DR = 0,
        \quad
        D' R^E = 0.
    \end{equation}
\end{lemma}
In the following we shall sometimes consider $\WL$ as sub-algebra of
$\WLEnd$. Hence $R \in \WLEnd$ satisfies $D'R = 0$, too. Moreover, for
an antisymmetric $k$-form $\Omega \in \Lambda^k \subseteq \WLO^k$ we
simply have
\begin{equation}
    \label{eq:DOmegas}
    D \Omega = D' \Omega = d\Omega.
\end{equation}

%
%

\section{The fiberwise deformations}
\label{sec:fiberdef}

For $f, g \in \CWeyl$ one defines the \emph{fiberwise Weyl product} by
\begin{equation}
    \label{eq:weylprod}
    f \circ g = \mu \circ \eu^{\frac{\im\lambda}{2}\Lambda^{kl}\,
      i_s(\partial_k) \otimes i_s(\partial_l)} f \otimes g,
\end{equation}
where $\mu$ is the undeformed fiberwise product and
$\Lambda^{kl} = - \omega^{kl}$ are the components of the Poisson
tensor. We extend this product to $\CWL$ and $\CWLEnd$ in the obvious
way. Finally we define
\begin{equation}
    \label{eq:leftmod}
    (f \otimes \alpha \otimes A) \circ (g \otimes \beta \otimes s)
    := f \circ g \otimes \alpha \wedge \beta \otimes As
\end{equation}
\begin{equation}
    \label{eq:rightmod}
    (g \otimes \beta \otimes s) \circ (h \otimes \gamma)
    := g \circ h \otimes \beta \wedge \gamma \otimes s
\end{equation}
for $f,g,h \in \CWeyl$, $\alpha, \beta \in \Lambda^k$, 
$A \in \Gamma^\infty(\End(E))$, and $s \in \Gamma^\infty(E)$ and
extend this by linearity to a left and right action of $\CWLEnd$ and
$\CWL$ on $\CWLE$, respectively. Then the following is obvious:
\begin{lemma}
    \label{lemma:fiberdeformations}
    The fiberwise Weyl product is globally well-defined and makes
    $\CWeyl$, $\CWL$ as well as $\CWLEnd$ into associative algebras
    with subalgebras $\Weyl$, $\WL$, and $\WLEnd$,
    respectively. Moreover, $\CWLE$ becomes a bimodule for $\CWLEnd$
    and $\CWL$ and $\WLE$ is a bimodule for $\WLEnd$ and $\WL$.
\end{lemma}
Moreover, $\circ$ is still (formally) $\dega$- and $\Deg$-graded,
i.e. $\dega$ and $\Deg$ are still derivations. The (super-) center
with respect to $\circ$ is now given by the anti-symmetric forms
$\Lambda^\bullet$. Moreover, a simple computation gives that
\begin{equation}
    \label{eq:deltaad}
    - \delta = \frac{\im}{\lambda} \ad (\womega)
\end{equation}
is an inner $\circ$-derivation. As we have chosen $\nabla$ to be
symplectic it turns out that $D$ as well as $D'$ are still derivations
of $\circ$. Moreover, $D^E$ is a module derivation in the sense that
\begin{equation}
    \label{eq:DEmodderivleft}
    D^E(a \circ \Psi) 
    = (D'a) \circ \Psi + (-1)^{\dega a} a \circ (D^E \Psi)
\end{equation}
\begin{equation}
    \label{eq:DEmodderivright}
    D^E(\Psi \circ b) 
    = (D^E \Psi) \circ b + (-1)^{\dega \Psi} \Psi \circ (Db)
\end{equation}
for all $a \in \CWLEnd$, $\Psi \in \CWLE$, and $b \in \CWL$ which is
the quantum analog of (\ref{eq:DDEDEndComp}). The squares of the
covariant derivatives turn out to be inner derivations with respect to
$\circ$:
\begin{lemma}
    \label{lemma:DsquaresInner}
    With $R$ and $R^E$ as in \eqref{eq:RDef} and \eqref{eq:REDef} we
    have
    \begin{equation}
        \label{eq:DsquareadR}
        D^2 = \frac{\im}{\lambda} \ad(R)
        \quad
        \textrm{and}
        \quad
        (D')^2 = \frac{\im}{\lambda} 
        \ad\left(R - \im \lambda R^E\right)
    \end{equation}
    \begin{equation}
        \label{eq:DEsquareInner}
        (D^E)^2 = \frac{\im}{\lambda} \ad(R) + R^E.
    \end{equation}
\end{lemma}
Here we note that $R$ can act from the left as well as from the right
on $\CWLE$ whence $\ad(R)$ is meaningful, while $R^E$ simply
acts by $\circ$-left multiplication. The proof for
\eqref{eq:DsquareadR} can be found in Fedosov's book
\cite[Sect.~5.3]{fedosov:1996} and \eqref{eq:DEsquareInner} is an easy
computation using Lemma~\ref{lemma:Dsquares}. Note that in
\eqref{eq:DsquareadR} the curvature $R^E$ appears with an additional
power of $\lambda$. This will play a major role later.

%
%

\section{The Fedosov derivatives}
\label{sec:fedder}

The main idea of Fedosov's construction is to realize
$C^\infty(M)[[\lambda]]$ and $\Gamma^\infty(\End(E))[[\lambda]]$ as
kernels of derivatives $\mathcal{D}$ and $\mathcal{D}'$ of the
fiberwise Weyl algebras $\WL$ and $\WLEnd$, respectively. The Ansatz
for $\mathcal{D}'$ (and analogously for $\mathcal{D}$) is
\begin{equation}
    \label{eq:calD}
    \mathcal{D}' = -\delta + D' + \frac{\im}{\lambda} \ad(r'),
\end{equation}
where $r' \in \CWLEnd$ is an element with total degree $\ge 3$ and
anti-symmetric degree $1$. Note that $\frac{\im}{\lambda}\ad(r')$
raises the $\Deg$-degree by at least $+1$ but it may lower the
$\lambda$-degree even if $r'$ does not contain negative powers of
$\lambda$ as in $\CWLEnd$ the undeformed product is \emph{not}
super-commutative. For the square of $\mathcal{D}'$ as in
\eqref{eq:calD} one has
\begin{equation}
    \label{eq:calDsquare}
    (\mathcal{D}')^2
    = \frac{\im}{\lambda} \ad 
    \left(
        -\omega - \delta r' 
        +  R - \im \lambda R^E
        + D'r' 
        + \frac{\im}{\lambda} r' \circ r' 
    \right),
\end{equation}
as a straightforward computation shows. Moreover, the `curvature' of
$\mathcal{D}'$ satisfies the Bianchi identity
\begin{equation}
    \label{eq:qbianchi}
    \mathcal{D}' 
    \left( 
        -\omega + \delta r' 
        + R - \im \lambda R^E
        + D' r' +
        \frac{\im}{\lambda} r' \circ r' 
    \right) = 0.
\end{equation}
Now if we want a `flat connection' $(\mathcal{D}')^2 = 0$ then the
curvature has to be a \emph{central} element, i.e. a formal power
series of two-forms 
$\Omega = \sum_{i=1}^\infty \lambda^i \Omega_i$. Then 
\eqref{eq:qbianchi} implies that necessarily $d\Omega = 0$. The
following theorem of Fedosov ensures that we can find such a $r'$ for
any given choice of $\Omega$.
\begin{theorem}[Fedosov\protect{\cite[Sect.~5.3]{fedosov:1996}}]
    \label{theorem:fedosovI}
    Let $\Omega = \sum_{i=1}^\infty \lambda^i \Omega_i$ be a closed
    two-form. Then there exists a unique $r' \in \CWLEnd$ with
    anti-symmetric degree $+1$ and total degree $\ge 3$ such that
    \begin{equation}
        \label{eq:rpime}
        \delta r' = 
        R - \im\lambda R^E + D' r' + \frac{\im}{\lambda} r' \circ r' 
        + \Omega
        \quad
        \textrm{and}
        \quad
        \delta^{-1} r' = 0. 
    \end{equation}
    In this case $(\mathcal{D}')^2 = 0$.
\end{theorem}
\begin{remark}
    \label{remark:neglpowers}
    A priori, the recursion for finding $r'$ works only in $\CWeyl$ and
    not in $\Weyl$ as the undeformed product of $\WLEnd$ is already
    non-commutative whence $\frac{\im}{\lambda} r' \circ r'$ may
    generate negative powers of $\lambda$. However, the crucial
    element causing the non-commutativity is $-\im\lambda R^E$ which
    comes with an additional power of $\lambda$. Then proof that $r'$
    does not contain negative powers of $\lambda$ can be done by
    induction using the following recursion formulas, see
    e.g.~\cite{bordemann.waldmann:1997a}.
\end{remark}
\begin{lemma}
    \label{lemma:rrecursion}
    The element $r'$ can be obtained recursively with respect to the
    total degree $\Deg$ by
    \begin{gather}
        \label{eq:rdrei}
        (r')^{(3)} = \delta^{-1} 
        \left(R - \im\lambda R^E + \lambda \Omega_1\right)
        \\
        \label{eq:rkplusdrei}
        (r')^{(k+3)} = \delta^{-1} 
        \left(
            D'(r')^{(k+2)} 
            + \frac{\im}{\lambda} \sum_{l=1}^{k-1}
            (r')^{(l+2)} \circ (r')^{(k+2-l)}
            + \left\{
                \begin{array}{cl}
                    \lambda^{k/2+1}\Omega_{k/2+1} & k \; \textrm{even} \\
                    0 & k \; \textrm{odd}
                \end{array}
            \right\}
        \right),
    \end{gather}
    where $r' = \sum_{k=3}^\infty (r')^{(k)}$. Moreover, 
    $r' \in \WLEnd$.
\end{lemma}

In a second step one computes the kernel of $\mathcal{D}'$. It turns
out that the kernel is in bijection to the sections 
$\Gamma^\infty (\End(E))[[\lambda]]$.
\begin{theorem}[Fedosov\protect{\cite[Sect.~5.3]{fedosov:1996}}]
    \label{theorem:fedosovII}
    The map 
    $\sigma: \ker \mathcal{D}' \cap \ker \dega \to
    \Gamma^\infty (\End(E))(\!(\lambda)\!)$ 
    is a $\mathbb{C}(\!(\lambda)\!)$-linear and $\lambda$-adically
    continuous bijection.
\end{theorem}
The inverse is denoted by $\tau'$ and referred to as the
\emph{Fedosov-Taylor series} as it is the quantum analog of the
formal Taylor series. For the Fedosov-Taylor series one has the
following recursion formula, see also~\cite{bordemann.waldmann:1997a}.
\begin{lemma}
    \label{lemma:taurecurs}
    The Fedosov-Taylor series of a section 
    $A \in \Gamma^\infty(\End(E))$ can be obtained recursively with
    respect to the total degree $\Deg$ by
    \begin{gather}
        \label{eq:taunull}
        \tau'(A)^{(0)} = A 
        \\
        \label{eq:taukpluseins}
        \tau'(A)^{(k+1)} = \delta^{-1} 
        \left(
            D' \tau'(A)^{(k)} + \frac{\im}{\lambda}
            \sum_{l=1}^{k-1} \ad\left((r')^{(l+2)}\right)
            \tau'(A)^{(k-l)} 
        \right),
    \end{gather}
    where $\tau'(A) = \sum_{k=0}^\infty \tau'(A)^{(k)}$. Moreover,
    $\tau(A) \in \WLEnd$.
\end{lemma}
From the two lemmas we observe that it is sufficient to stay within the
framework of formal power series in $\lambda$: we do not need the
extension to $\CWLEnd$ \emph{a posteriori}. Note however, that in the
original recursion it is not so obvious that we do not produce
negative $\lambda$-powers.

Since the kernel of a super-derivation is a sub-algebra one can
pull-back the fiberwise Weyl product $\circ$ of $\WLEnd$ to
$\Gamma^\infty(\End(E))[[\lambda]]$ by means of $\sigma$ and
$\tau'$. Hence we obtain an associative deformation
\begin{equation}
    \label{eq:starendos}
    A \star' B = \sigma \left(\tau'(A) \circ \tau'(B)\right)
\end{equation}
for $A, B \in \Gamma^\infty(\End(E))[[\lambda]]$.

Of course the same line of argument can be applied to $\WL$ itself
without endomorphism-valued elements. In this case we obtain a unique
$r \in \WLO^1$ of total degree $\ge 3$ such that
\begin{equation}
    \label{eq:deltar}
    \delta r = R + Dr + \frac{\im}{\lambda} r \circ r + \Omega
    \quad
    \textrm{and}
    \quad
    \delta^{-1}r = 0
\end{equation}
with corresponding Fedosov derivative
$\mathcal{D} = -\delta + D + \frac{\im}{\lambda} \ad(r)$ and
$\mathcal{D}^2 = 0$. Note that we have used the \emph{same} closed
two-form $\Omega$ as in \eqref{eq:rpime} to specify $r$ and thus
$\mathcal{D}$. Again one has a Fedosov-Taylor series 
\begin{equation}
    \label{eq:fedtaylor}
    \tau: C^\infty (M)[[\lambda]] \to \ker\mathcal{D} \cap \ker\dega
    \quad
    \textrm{with}
    \quad
    \sigma(\tau(f)) = f
\end{equation}
and a corresponding associative deformation
\begin{equation}
    \label{eq:starprod}
    f \star g = \sigma(\tau(f) \circ \tau(g)).
\end{equation}
It turns out that this is actually a star product. Up to now we have
just recalled Fedosov's original construction of the star product
$\star$ and the associative deformation $\star'$ of 
$\Gamma^\infty(\End(E))$.

For later use we shall consider the two elements $r$ and $r'$ more
closely:
\begin{lemma}
    \label{lemma:rE}
    The classical limits of $r$ and $r'$ coincide whence
    \begin{equation}
        \label{eq:rEDef}
        r^E := \frac{\im}{\lambda} (r' - r) \in \WLEnd
    \end{equation}
    does not contain negative powers of $\lambda$. Moreover, $r^E$ is
    uniquely determined by 
    \begin{equation}
        \label{eq:rEdelta}
        \delta r^E = R^E + D'r^E + \frac{\im}{\lambda} \ad(r) 
        r^E + r^E \circ r^E
        \quad
        \textrm{and}
        \quad
        \delta^{-1} r^E = 0.
    \end{equation}
    In particular, $r^E = 0$ if and only if $\nablaE$ is flat,
    i.e. $R^E = 0$.
\end{lemma}
\begin{proof}
    First note that $D'r = Dr$. Comparing the recursion formulas for
    $r$ and $r'$ we see that their difference has to be of order
    $\lambda$. Moreover, \eqref{eq:rEdelta} follows directly from
    \eqref{eq:rpime} and \eqref{eq:deltar} since we have used the
    \emph{same} $\Omega$. Finally, any element $a$ of anti-symmetric
    degree $+1$ is uniquely determined by specifying $\delta a$ and
    $\delta^{-1}a$ according to \eqref{eq:hodge}. Clearly $r^E = 0$
    implies $R^E = 0$ as this is the component of $r^E$ of total
    degree $1$. Conversely, if $R^E = 0$ then $r^E = 0$ follows as in
    this case it is the unique solution of \eqref{eq:rEdelta}.
\end{proof}

%
%

\section{The bimodule structure and Morita equivalence}
\label{sec:bimodule}

From the construction of $\star$ and $\star'$ it is easy to guess how
one can deform the classical bimodule structure of
$\Gamma^\infty(E)$. We just have to find a suitable Fedosov derivative
$\mathcal{D}^E$ with Fedosov-Taylor series $\tau^E$ for $\WLE$,
too. This can indeed be done, we even have already all pieces present
and do not have to start a new recursion. We define 
$\mathcal{D}^E: \WLE \to \WLO^{\!\bullet+1}\!\otimes\!\mathcal{E}$
by
\begin{equation}
    \label{eq:calDEDef}
    \mathcal{D}^E = - \delta + D^E + \frac{\im}{\lambda} \ad(r) + r^E,
\end{equation}
where $r$ and $r^E$ are given as before. Note that $\ad(r)$ is
well-defined as $r$ can act form left and right. The element $r^E$ is
understood to act by $\circ$-left multiplication. Finally note, that
$\ad(r)\Psi$ is always of order $\lambda$, whence $\mathcal{D}^E$ does
not produce negative powers of $\lambda$.
\begin{theorem}
    \label{theorem:DE}
    The Fedosov derivative $\mathcal{D}^E$ satisfies
    \begin{gather}
        \label{eq:DEModDerifleft}
        \mathcal{D}^E(a \circ \Psi) 
        = (\mathcal{D}' a) \circ \Psi 
        + (-1)^{\dega a} a \circ (\mathcal{D}^E\Psi) 
        \\
        \label{eq:DEModDerifright}
        \mathcal{D}^E(\Psi \circ b) 
        = (\mathcal{D}^E \Psi) \circ b
        + (-1)^{\dega \Psi} \Psi \circ (\mathcal{D}b)      
    \end{gather}
    as well as $(\mathcal{D}^E)^2 = 0$. Moreover,
    \begin{equation}
        \label{eq:sigmaE}
        \sigma: \ker\mathcal{D}^E \cap \ker\dega 
        \to \Gamma^\infty(E)[[\lambda]]
    \end{equation}
    is a $\mathbb{C}[[\lambda]]$-linear bijection with inverse denoted
    by $\tau^E$. The Fedosov-Taylor series $\tau^E(s)$ of a section 
    $s \in \Gamma^\infty(E)$ can be obtained recursively with respect
    to the total degree by
    \begin{gather}
        \label{eq:tauErecursNull}
        \tau^E(s)^{(0)} = s \\
        \tau^E(s)^{(k+1)} = \delta^{-1}
        \left(D^E \tau^E(s)^{(k)} +
            \sum_{l=1}^k
            \left(
                \frac{\im}{\lambda}\ad\left(r^{(l+2)}\right) +
                (r^E)^{(l)} 
            \right)
            \tau^E(s)^{(k-l)}
        \right)
    \end{gather}
    where $\tau^E(s) = \sum_{k=0}^\infty \tau^E(s)^{(k)}$.
\end{theorem}
\begin{proof}
    The equations \eqref{eq:DEModDerifleft},
    \eqref{eq:DEModDerifright} and $(\mathcal{D}^E)^2 = 0$ are just
    straightforward computations using the relations between $r$ and
    $r^E$ as well as the results from Section~\ref{sec:fiberdef}. The
    crucial point is that we have used the \emph{same} $\Omega$ for
    $r'$ and $r$. Now let $s \in \Gamma^\infty(E)[[\lambda]]$ be given
    and define the operator
    \begin{equation}
        \label{eq:Ts}
        T_s \Psi := s + \delta^{-1} 
        \left(D^E \Psi + \frac{\im}{\lambda} \ad(r)\Psi 
            + r^E \circ \Psi
        \right)
    \end{equation}
    for $\Psi \in \WLE$ of anti-symmetric degree $0$. It immediately
    follows that $T_s$ is strictly contracting in the complete metric
    space $\WLO^{\!0}\!\otimes\!\mathcal{E}$, where the ultra-metric is
    defined by means of the $\Deg$-degree. Thus $T_s$ has a unique
    fixed point, denoted by $\tau^E(s)$, by Banach's fixed point
    theorem, see
    e.g.~\cite[App.~A]{bordemann.neumaier.waldmann:1999}. One also
    obtains that $s \mapsto \tau^E(s)$ is
    $\mathbb{C}[[\lambda]]$-linear. From \eqref{eq:Ts} it follows that
    $\delta^{-1} \mathcal{D}^E\tau^E(s) = 0$ and with
    $(\mathcal{D}^E)^2 = 0$ it follows that
    $\delta\mathcal{D}^E\tau^E(s) = (D^E + \frac{\im}{\lambda}\ad(r) +
    r^E)\tau^E(s)$ whence by \eqref{eq:hodge} we have
    \begin{equation}
        \label{eq:DeTauE}
        \mathcal{D}^E \tau^E(s) 
        = \delta^{-1} 
        \left( D^E + \frac{\im}{\lambda} \ad(r) + r^E\right) 
        \mathcal{D}^E \tau^E(s).
    \end{equation}
    Thus $\mathcal{D}^E\tau^E(s)$ turns out to be the fixed point of a
    strictly contracting \emph{linear} operator, whence
    $\mathcal{D}^E\tau^E(s) = 0$. On the other hand, let
    $\mathcal{D}^E\Psi = 0$. Then with \eqref{eq:hodge} it follows
    that $\Psi$ is the fixed point of $T_{\sigma(\Psi)}$ whence 
    $\Psi = \tau^E(\sigma(\Psi))$. Finally, the recursion formulas
    follow immediately from $\tau^E(s) = T_s \tau^E(s)$.
\end{proof}
\begin{corollary}
    \label{corollary:bimodule}
    The sections $\Gamma^\infty(E)[[\lambda]]$ become a bimodule for
    $\star'$ and $\star$ by
    \begin{equation}
        \label{eq:bimodDef}
        A \bullet' s = \sigma (\tau'(A) \circ \tau^E(s))
        \quad
        \textrm{and}
        \quad
        s \bullet f = \sigma (\tau^E(s) \circ \tau(f)),
    \end{equation}
    where $A \in \Gamma^\infty(\End(E))[[\lambda]]$,
    $s \in \Gamma^\infty(E)[[\lambda]]$, and
    $f \in C^\infty(M)[[\lambda]]$.
\end{corollary}

Let us now focus on the case of a line bundle $E = L \to M$. In this
case $\Gamma^\infty(\End(L))[[\lambda]] = C^\infty(M)[[\lambda]]$
whence the product $\star'$ is defined for \emph{functions} on $M$. In
fact, it turns out to be a star product as well.

By definition, the \emph{Fedosov class} of a Fedosov star product is
the deRham cohomology class of the curvature of the corresponding
Fedosov derivative, i.e.
\begin{equation}
    \label{eq:Fclass}
    F(\star) = 
    \left[
        \omega + \delta r - R - Dr - \frac{\im}{\lambda} r \circ r
    \right]
    = [\omega] + [\Omega],
\end{equation}
if $\star$ is obtained from 
$\mathcal{D} = - \delta + D + \frac{\im}{\lambda} \ad(r)$. Thus,
according to \eqref{eq:rpime}, the class of $\star'$ is given by
\begin{equation}
    \label{eq:Fclassprime}
    F(\star') =  
    \left[
        \omega + \delta r' - R - Dr' - \frac{\im}{\lambda} r' \circ r'
    \right]
    = [\omega] + [\Omega] + \im\lambda[R^L],
\end{equation}
where $R^L$ is the curvature of $\nabla^{\scriptscriptstyle L}$. But
this is just the \emph{Chern class} of $L$ whence
\begin{equation}
    \label{eq:FclassL}
    F(\star') =  F(\star) - 2 \pi\lambda c_1(L).
\end{equation}
Taking into account that the \emph{characteristic class} $c(\star)$ of
a Fedosov star product is given by $\frac{1}{\im\lambda} F(\star)$,
see e.g. the discussion in \cite{neumaier:1999:pre}, we have the
following corollary:
\begin{corollary}
    \label{corollary:MEclass}
    In case of a line bundle $E = L \to M$, the characteristic classes
    of $\star'$ and $\star$ are related by
    \begin{equation}
        \label{eq:MEclass}
        c(\star') = c(\star) + 2\pi\im c_1(L).
    \end{equation}
\end{corollary}
\begin{remark}
    \label{remark:ME}
    This is of course to be expected from
    \cite[Thm.~3.1]{bursztyn.waldmann:2001:pre} as the bimodule
    structure $\bullet'$ and $\bullet$ on
    $\Gamma^\infty(L)[[\lambda]]$ is exactly a Morita equivalence
    bimodule for the two star products $\star'$ and $\star$. The
    remarkable point is that the computation of $c(\star')$ is almost
    a triviality in the Fedosov framework, compared to the \v{C}ech
    cohomological computation in \cite{bursztyn.waldmann:2001:pre}.
\end{remark}

%
%

\section{The case of a Hermitian vector bundle}
\label{sec:herm}

Consider now a Hermitian fiber metric $h$ for $E$ and assume that
$\nablaE$ is compatible with $h$. In this case
$\Gamma^\infty(\End(E))$ has a natural $^*$-involution
defined by $h(As, s') = h(s, A^*s')$ for 
$A \in \Gamma^\infty(\End(E))$ and
$s, s' \in \Gamma^\infty(E)$. Thus we can extend this $^*$-involution,
together with the complex conjugation, to a super-$^*$-involution of
$\WLEnd$ and $\WL$, respectively. It is well-known that the fiberwise
Weyl-product $\circ$ is compatible with this $^*$-involution, 
i.e. $(a \circ b)^* = (-1)^{\dega a \, \dega b} \, b^* \circ a^*$ for
all $a,b \in \WLEnd$. By the unique characterization of $r'$ and $r$
by \eqref{eq:rpime} and \eqref{eq:deltar}, respectively, the following
lemma is straightforward, see also 
\cite[Lem.~3.3]{bordemann.waldmann:1997a}.
\begin{lemma}
    \label{lemma:Reality}
    Let $\Omega = \cc{\Omega}$ be a real formal two-form. Then
    \begin{equation}
        \label{eq:realrrr}
        (r')^* = r',
        \quad
        \cc{r} = r,
        \quad
        (r^E)^* = - r^E,
    \end{equation}
    as well as
    \begin{equation}
        \label{eq:realDD}
        (\mathcal{D}'a)^* = \mathcal{D}' a^*,
        \quad
        \textrm{and}
        \quad
        \cc{(\mathcal{D} b)} = \mathcal{D} \cc{b}
    \end{equation}
    for $a \in \WLEnd$ and $b \in \WL$. Moreover,
    \begin{equation}
        \label{eq:realtau}
        \tau'(A^*) = (\tau'(A))^*,
        \quad
        \tau(\cc{f}) = \cc{\tau(f)},
    \end{equation}
    whence $\star'$ and $\star$ are Hermitian deformations, i.e.
    \begin{equation}
        \label{eq:starHerm}
        (A \star' B)^* = B^* \star' A^*
        \quad
        \textrm{and}
        \quad
        \cc{f \star g} = \cc{g} \star \cc{f}.
    \end{equation}
\end{lemma}
Let us assume $\Omega = \cc{\Omega}$ for the following. Then in a next
step we extend the fiber metric to $\WLE$ with values in $\WL$ by
defining
\begin{equation}
    \label{eq:fiberh}
    H (f \otimes \alpha \otimes s, g \otimes \beta \otimes s')
    =
    \cc{f} \circ g \otimes \cc{\alpha} \wedge \beta h(s, s')
\end{equation}
and extending this by sesquilinearity. Note that we can write
$h(s,s')$ on any side of the tensor product or of the $\circ$-product
as it is a \emph{function}. The following properties are immediate:
\begin{lemma}
    \label{lemma:H}
    Let $a \in \WLEnd$, $\Psi,\Psi' \in \WLE$, and $b \in \WL$. Then
    \begin{equation}
        \label{eq:Hlinearity}
        H(a \circ \Psi, \Psi') = H(\Psi, a^* \circ \Psi')
        \quad
        \textrm{and}
        \quad
        H(\Psi, \Psi' \circ b) = H(\Psi, \Psi') \circ b
    \end{equation}
    as well as
    \begin{equation}
        \label{eq:HHerm}
        H(\Psi, \Psi') = \cc{H(\Psi', \Psi)}.
    \end{equation}
\end{lemma}
Using Lemma~\ref{lemma:Reality} and \ref{lemma:H} as well as the
properties of $\mathcal{D}^E$ as in \eqref{eq:DEModDerifleft} and
\eqref{eq:DEModDerifright} we obtain by a simple computation the
following compatibility
\begin{equation}
    \label{eq:DHcomp}
    \mathcal{D} (H(\Psi, \Psi')) = 
    H (\mathcal{D}^E \Psi, \Psi') 
    + (-1)^{\dega \Psi} H(\Psi, \mathcal{D}^E \Psi')
\end{equation}
where $\Psi, \Psi' \in \WLE$. Clearly this can be seen as a direct
analog of the compatibility of $\nablaE$ and $h$. As a consequence
we can define a deformed Hermitian metric $\boldsymbol{h}$ by
\begin{equation}
    \label{eq:defh}
    \boldsymbol{h} (s, s') = \sigma(H(\tau^E(s), \tau^E(s'))).
\end{equation}
\begin{theorem}
    \label{theorem:Defh}
    The map $\boldsymbol{h}$ is $\mathbb{C}[[\lambda]]$-sesquilinear
    and satisfies
    \begin{gather}
        \label{eq:defhHerm}
        \boldsymbol{h}(s, s') = \cc{\boldsymbol{h}(s',s)}
        \quad
        \textrm{and}
        \quad
        \boldsymbol{h}(s, s) \ge 0
        \\
        \label{eq:defhlinear}
        \boldsymbol{h} (s, s'\bullet f) = 
        \boldsymbol{h}(s, s') \star f
        \\
        \label{eq:defhAdjoints}
        \boldsymbol{h} (A \bullet' s, s')
        = \boldsymbol{h}(s, A^* \bullet' s')
    \end{gather}
    for all $s, s' \in \Gamma^\infty(E)[[\lambda]]$, 
    $f \in C^\infty(M)[[\lambda]]$ and $A \in
    \Gamma^\infty(\End(E))[[\lambda]]$.
\end{theorem}
\begin{proof}
    The only non-trivial statement is the positivity 
    $\boldsymbol{h}(s, s) \ge 0$ which follows from
    \cite[Prop.~2.8]{bursztyn.waldmann:2000b}. The other statements
    are straightforward computations using (\ref{eq:Hlinearity}),
    (\ref{eq:HHerm}), and (\ref{eq:DHcomp}).
\end{proof}

%
%

\section{Conclusion and further questions}
\label{sec:conclusion}

Let us conclude with a few remarks on this approach to Morita
equivalence of Fedosov star products. First we would like to point out
that Corollary~\ref{corollary:MEclass}, together with the fact that
every star product is equivalent to some Fedosov star product, implies
that the condition \eqref{eq:MEclass} is actually sufficient for
Morita equivalence. To prove that \eqref{eq:MEclass} is also
necessary, one has to go beyond Fedosov's construction as within this
construction it is not directly clear that any bimodule deformation is
equivalent to the one we described here. This fact was shown in
\cite{bursztyn.waldmann:2000b,bursztyn.waldmann:2001:pre}.

Nevertheless, the above construction is very geometric and almost
explicit which makes it attractive to consider the following
questions:
\begin{enumerate}
\item The whole construction depends functorially on $\nabla$ and
    $\nablaE$, say for the choice $\Omega = 0$. Thus having a group
    acting on $M$ with a lift to an action on the vector bundle $E$
    such that $\nabla$ and $\nablaE$ are preserved, this will lead to
    invariant $\star$, $\star'$, $\bullet$, $\bullet'$. A more
    detailed investigation of this problem will be subject to a future
    project.
\item On a K{\"a}hler manifold one has a Fedosov construction for a
    star product of Wick type (separation of variables),
    see~\cite{karabegov:1996,bordemann.waldmann:1997a,karabegov:2000}. 
    Here one should be able to find a similar notion of `separation of
    variables' for $\bullet$ and $\bullet'$ for the case of
    holomorphic vector bundles. Moreover, even in the case of almost
    K{\"a}hler manifolds \cite{karabegov.schlichenmaier:2001}, one can
    try to understand the characteristic classes of the star products
    of Wick type as coming from a particular bimodule deformation.
    Again this will be discussed in a forthcoming project.
\item Finally, we hope that this construction can be transfered to the
    case of Poisson manifolds using the `Fedosov-like' approach of
    Cattaneo, Felder and Tomassini
    \cite{cattaneo.felder.tomassini:2001:pre} in order to give more
    insight in the classification of star products up to Morita
    equivalence on Poisson manifolds. It may also give a alternative,
    global description of the appearance of Morita equivalence in
    non-commutative field theories as discussed in
    \cite{jurco.schupp.wess:2001b:pre}.
\end{enumerate}

%
%


\begin{thebibliography}{10}

\bibitem {bayen.et.al:1978}
{\sc Bayen, F., Flato, M., Fr{{\o}}nsdal, C., Lichnerowicz, A., Sternheimer,
  D.: }\newblock {\em Deformation Theory and Quantization}.
\newblock Ann. Phys.  {\bf 111} (1978), 61--151.

\bibitem {bertelson.cahen.gutt:1997}
{\sc Bertelson, M., Cahen, M., Gutt, S.: }\newblock {\em Equivalence of Star
  Products}.
\newblock Class. Quantum Grav.  {\bf 14} (1997), A93--A107.

\bibitem {bordemann.neumaier.waldmann:1999}
{\sc Bordemann, M., Neumaier, N., Waldmann, S.: }\newblock {\em Homogeneous
  Fedosov star products on cotangent bundles II: GNS representations, the WKB
  expansion, traces, and applications}.
\newblock J. Geom. Phys.  {\bf 29} (1999), 199--234.

\bibitem {bordemann.waldmann:1997a}
{\sc Bordemann, M., Waldmann, S.: }\newblock {\em A Fedosov Star Product of
  Wick Type for K{\"{a}}hler Manifolds}.
\newblock Lett. Math. Phys.  {\bf 41} (1997), 243--253.

\bibitem {bordemann.waldmann:1998}
{\sc Bordemann, M., Waldmann, S.: }\newblock {\em Formal GNS Construction and
  States in Deformation Quantization}.
\newblock Commun. Math. Phys.  {\bf 195} (1998), 549--583.

\bibitem {bursztyn.waldmann:2000b}
{\sc Bursztyn, H., Waldmann, S.: }\newblock {\em Deformation Quantization of
  Hermitian Vector Bundles}.
\newblock Lett. Math. Phys.  {\bf 53} (2000), 349--365.

\bibitem {bursztyn.waldmann:2001a}
{\sc Bursztyn, H., Waldmann, S.: }\newblock {\em Algebraic Rieffel Induction,
  Formal Morita Equivalence and Applications to Deformation Quantization}.
\newblock J. Geom. Phys.  {\bf 37} (2001), 307--364.

\bibitem {bursztyn.waldmann:2001:pre}
{\sc Bursztyn, H., Waldmann, S.: }\newblock {\em The characterisitc classes of
  {M}orita equivalence of star products on symplectic manifolds}.
\newblock Preprint Freiburg FR-THEP 2001/09,  {\bf math.QA/0106178} (June
  2001).

\bibitem {cattaneo.felder.tomassini:2001:pre}
{\sc Cattaneo, A.~S., Felder, G., Tomassini, L.: }\newblock {\em Fedosov
  connections on jet bundles and deformation quantization}.
\newblock Preprint  {\bf math.QA/0111290} (November 2001).

\bibitem {dewilde.lecomte:1983b}
{\sc DeWilde, M., Lecomte, P. B.~A.: }\newblock {\em Existence of Star-Products
  and of Formal Deformations of the Poisson Lie Algebra of Arbitrary Symplectic
  Manifolds}.
\newblock Lett. Math. Phys.  {\bf 7} (1983), 487--496.

\bibitem {dito.sternheimer:2001:pre}
{\sc Dito, G., Sternheimer, D.: }\newblock {\em Deformation quantization:
  genesis, developments and metamorphoses}.
\newblock Preprint   (2001).
\newblock Contribution to the Proceedings for the {$68^{\textrm{eme}}$}
  Rencontre entre Physiciens Theoriciens et Mathematiciens on Deformation
  Quantization. Strasbourg, 31.~05.~2001 -- 02.~06.~2001.

\bibitem {fedosov:1996}
{\sc Fedosov, B.~V.: }\newblock {\em Deformation Quantization and Index
  Theory}.
\newblock Akademie Verlag, Berlin, 1996.

\bibitem {gutt:2000}
{\sc Gutt, S.: }\newblock {\em Variations on deformation quantization}.
\newblock In: {\sc Dito, G., Sternheimer, D. (eds.): }\newblock {\em
  Conf{\`e}rence Mosh{\`e} Flato 1999. Quantization, Deformations, and
  Symmetries}, {\em Mathematical Physics Studies} no. {\bf 21},   217--254.
  Kluwer Academic Publishers, Dordrecht, Boston, London, 2000.

\bibitem {jurco.schupp.wess:2001b:pre}
{\sc Jurco, B., Schupp, P., Wess, J.: }\newblock {\em Noncommutative line
  bundle and {M}orita equivalence}.
\newblock Preprint  {\bf hep-th/0106110} (June 2001).

\bibitem {karabegov:1996}
{\sc Karabegov, A.~V.: }\newblock {\em Deformation Quantization with Separation
  of Variables on a K{\"{a}}hler Manifold}.
\newblock Commun. Math. Phys.  {\bf 180} (1996), 745--755.

\bibitem {karabegov:2000}
{\sc Karabegov, A.~V.: }\newblock {\em On Fedosov's approach to Deformation
  Quantization with Separation of Variables}.
\newblock In: {\sc Dito, G., Sternheimer, D. (eds.): }\newblock {\em
  Conf{\`e}rence Mosh{\`e} Flato 1999. Quantization, Deformations, and
  Symmetries}, {\em Mathematical Physics Studies} no. {\bf 22}. Kluwer Academic
  Publishers, Dordrecht, Boston, London, 2000.

\bibitem {karabegov.schlichenmaier:2001}
{\sc Karabegov, A.~V., Schlichenmaier, M.: }\newblock {\em
  Almost-{K}{{\"a}}hler Deformation Quantization}.
\newblock Lett. Math. Phys.  {\bf 57} (2001), 135--148.

\bibitem {kontsevich:1997:pre}
{\sc Kontsevich, M.: }\newblock {\em Deformation Quantization of Poisson
  Manifolds, I}.
\newblock Preprint  {\bf q-alg/9709040} (September 1997).

\bibitem {nest.tsygan:1995a}
{\sc Nest, R., Tsygan, B.: }\newblock {\em Algebraic Index Theorem}.
\newblock Commun. Math. Phys.  {\bf 172} (1995), 223--262.

\bibitem {neumaier:1999:pre}
{\sc Neumaier, N.: }\newblock {\em Local $\nu$-Euler Derivations and Deligne's
  Characteristic Class of Fedosov Star Products}.
\newblock Preprint Freiburg FR-THEP-99/3,  {\bf math.QA/9905176} (May 1999).

\bibitem {omori.maeda.yoshioka:1991}
{\sc Omori, H., Maeda, Y., Yoshioka, A.: }\newblock {\em Weyl Manifolds and
  Deformation Quantization}.
\newblock Adv. Math.  {\bf 85} (1991), 224--255.

\bibitem {rieffel:1974b}
{\sc Rieffel, M.~A.: }\newblock {\em Morita equivalence for {$C^*$}-algebras
  and {$W^*$}-algebras}.
\newblock J. Pure. Appl. Math.  {\bf 5} (1974), 51--96.

\bibitem {waldmann:2001b:pre}
{\sc Waldmann, S.: }\newblock {\em On the Representation Theory of Deformation
  Quantization}.
\newblock Preprint Freiburg FR-THEP 2001/10,  {\bf math.QA/0107112} (July
  2001).
\newblock Contribution to the Proceedings for the {$68^{\textrm{eme}}$}
  Rencontre entre Physiciens Theoriciens et Mathematiciens on Deformation
  Quantization. Strasbourg, 31.~05.~2001 -- 02.~06.~2001.

\bibitem {weinstein:1994}
{\sc Weinstein, A.: }\newblock {\em Deformation Quantization}.
\newblock S\'eminaire Bourbaki 46\`eme ann\'ee  {\bf 789} (1994).

\end{thebibliography}
\begin{small}

\end{small}
\end{document}